\documentclass[a4paper,11pt]{revtex4-1}
\usepackage{amsmath}
\usepackage{amsfonts}
\usepackage{amssymb}
\usepackage{graphicx}
\usepackage{here}

\begin{document}

\title{Functional Dynamics by Intention Recognition in Iterated Games}
\author{Yuma Fujimoto$^{\dag}$ and Kunihiko Kaneko}
\affiliation{Department of Basic Science, The University of Tokyo, 3-8-1 Komaba, Meguro-ku, Tokyo 153-8902, Japan}
\affiliation{$^{\dag}${\rm yfujimoto@complex.c.u-tokyo.ac.jp}}
\date{\today}

\begin{abstract}
{Intention recognition is an important characteristic of intelligent agents. In their interactions with others, they try to read others' intentions and make an image of others to choose their actions accordingly. While the way in which players choose their actions depending on such intentions has been investigated in game theory, how dynamic changes in intentions by mutually reading others' intentions are incorporated into game theory has not been explored.

We present a novel formulation of game theory in which players read others' intentions and change their own through an iterated game. Here, intention is given as a function of the other's action and the own action to be taken accordingly as the dependent variable, while the mutual recognition of intention is represented as the functional dynamics.

It is shown that a player suffers no disadvantage when he/she recognizes the other's intention, whereas the functional dynamics reach equilibria in which both players' intentions are optimized. These cover a classical Nash and Stackelberg equilibria but we extend them in this study: Novel equilibria exist depending on the degree of mutual recognition. Moreover, the degree to which each player recognizes the other can also differ. This formulation is applied to resource competition, duopoly, and prisoner's dilemma games. For example, in the resource competition game with player-dependent capacity on gaining the resource, the superior player's recognition leads to the exploitation of the other, while the inferior player's recognition leads to cooperation through which both players' payoffs increase.}
\end{abstract}

\maketitle


\section{Introduction}
How each individual decides his/her own behavior is a long-standing problem in nature. Each agent freely behaves and receives a reward as a result. This situation is generally formulated as a ``normal form game'' in which the ``player,'' ``action,'' and ``payoff'' are given \cite{Neumann1944}. In a standard game played only once, each player needs to decide what action he/she chooses as his/her own strategy. An optimal strategy is to choose the best rewarded action depending on others' actions, which results in the Nash equilibrium \cite{Nash1950, Nash1951}. However, while this is the optimal strategy in a one-shot game, it is common that the game is repeatedly played in reality. In a repeated game, each player can refer to a data set on the history of actions played in the past and use this for the next game or a later one. The Nash equilibrium was not originally introduced to deal with such a situation. Instead, one possible equilibrium is given by the ``folk theorem'' \cite{Aumann1994, Rubinstein1980, Fudenberg2009}, where deviation from the Nash equilibrium can happen; however, the theorem only provides the requirements for the achieved equilibrium and cannot specify which equilibrium is really achieved. To determine an optimal strategy, which is now nothing given as in a one-shot game, we must therefore adopt a concrete learning process through which each player improves his/her own strategy against that of the other player.

An important characteristic of an intelligent agent (i.e., a human) is to recognize and make an image of others by using the history of the other's actions. Such an agent assumes that the other intentionally changes the next action in response to the agent's own action. For example, a descriptive and predictive model for a person's cognitive behavior has been proposed \cite{Nagel1995, Camerer2004}, as is based on the experiments of a repeated beauty contest game \cite{Keynes1936, Moulin1986}. In this model, each person is given a cognitive level. The level 0 person chooses an action with no recognition (randomly), while the level $k(>0)$ person best responds to the image of others at the $(k-1)$th level (or lower). Recently, such cognitive levels have been uncovered in neural economics \cite{Coricelli2009}. As another example, it has been reported that information on others' true intentions increases the performance in the game when the actions of others are transferred with disturbances by noise \cite{Charness2007, Falk2008, Fischbacher2013, Rand2015}. Following this experimental result, the evolution of the ability to recognize others' intentions is theoretically discussed \cite{Anh2011, Anh2015, Nakamura2016}. Although the above studies sufficiently justify the existence of humans' ability to recognize others and benefit from that, a theory for the dynamic coevolution of images between agents has been underdeveloped. In \cite{Taiji1999}, how a player builds and deconstructs the other's image in prisoner's dilemma games has been studied by using a recurrent network. How the equilibrium of actions is shaped and how it deviates from the Nash equilibrium, however, are not analyzed. Here, we develop a theoretical framework with mutual learning that shapes the other's internal intentions, generally applicable to any games, without resorting to specific learning algorithms.

Once every agent has an image of the other and best responds to it, the sequential actions are given. Therefore, how an agent constructs the other's image itself is now a strategy in repeated games, which is represented as a function, as shown later. Initially, both agents best respond to each other without referring to the other's strategy. Hence, both best-response functions as strategies achieve nothing but the Nash equilibrium \cite{Nash1950, Nash1951}. Then, before considering the dynamics of a pair of such strategy functions, we see an extreme case that an agent one-sidedly reads the other's strategy function. In this case, the Stackelberg equilibrium \cite{Stackelberg2011} is achieved, which is defined as an equilibrium in an ``extensive form game'' \cite{Kuhn1950} under perfect information.

Following these introductory results, we study the dynamics of strategy functions, which represent that agents mutually recognize the other's intention. With repeated games, each agent accurately reads the other's strategy function and optimizes his/her own one based on it. This dynamics reaches an equilibrium when there is no additional advantage for the further recognition of the other's strategy. At this point, a ``functional equilibrium'' is achieved between both players' strategies instead of the original Nash equilibrium.

Note that our formulation can be applied to general games. Here, applications to resource competition, duopoly, and probabilistic prisoner's dilemma games are provided as examples. In the former case, it is found that learning by an inferior agent increases the payoffs of both players, while that by a superior agent enhances exploitation and decreases the payoff of the other.

\section{Nash equilibrium}
We consider a two-player game in which players are denoted by $i\in\{1,2\}$. In addition, each player $i$'s action and its payoff are represented by $x_i$ and $u_i(x_1,x_2)$, respectively, which are continuous variables. A player tries to receive a higher payoff by optimizing his/her action depending on the other's action. Now, each player has an intention on which action he/she chooses depending on the other's action. The intentions vary depending on how the player imagines the other's action. Thus, the intention of player 1 is given as strategy function $f_1(x_2)$, which represents that action $x_1$ is chosen when player 2 takes action $x_2$. Player 2's strategy function is similarly defined as $f_2(x_1)$. Then, assuming that each player's action follows his/her own strategy, the equilibrium set of actions, denoted by $(x_1^{\mathrm{eq}},x_2^{\mathrm{eq}})$, is given by the crossing point of both players' functions. In other words, we get
\begin{equation}
\begin{split}
	&x_1^{\mathrm{eq}}=f_1(x_2^{\mathrm{eq}}),\\
	&x_2^{\mathrm{eq}}=f_2(x_1^{\mathrm{eq}}).
\end{split}
\label{eq1-00}
\end{equation}

In this section, we consider a situation in which both players have no recognition of the other's intention. In this case, each player simply maximizes his/her own payoff without referring to the other's strategy. To be consistent with the standard terminology, this strategy is called the ``best-response'' \cite{Nash1951}, as denoted by $f_1^{\mathrm{B}}(x_2)$ (B-response) for player 1. According to this definition, $f_1^{\mathrm{B}}(x_2)$ satisfies
\begin{equation}
	f_1^{\mathrm{B}}(x_2):=\mathrm{argmax}_{x_1}u_1(x_1,x_2).
\label{eq1-01}
\end{equation}

Eq.~\ref{eq1-01} simply means that player 1's strategy function is given by maximizing the payoff under the assumption that the other's action is constant independent of his/her own action. At this point, note that the strategy function given by B-response depends on the other's action $x_2$.

Player 2's B-response is given in the same way. Thus, when both players make B-responses, the equilibrium, denoted by $(x_1^{\mathrm{BB}},x_2^{\mathrm{BB}})$, is nothing but the Nash equilibrium from its definition. In the present paper, however, we call it the BB equilibrium, where the left index indicates the player's strategy to the other's strategy given by the right index, because the same equilibrium set of actions can be achieved by different sets of functions. In this study, which pair of functions results in the equilibrium action is important; hence, we need to specify not only the equilibrium point but the pair of functions to achieve it. At the BB equilibrium, each player's payoff is defined as $u_i^{\mathrm{BB}}:=u_i(x_1^{\mathrm{BB}},x_2^{\mathrm{BB}})$.

\section{Definition of the learning response and one-sided recognition}
Next, we define another type of intention where a player perfectly recognizes the other's intention. Then, each player optimizes his/her action based on the information on the other's strategy function. This strategy is termed the ``learning response'' (L-response), denoted by $f_1^{\mathrm{L}}(x_2)$, which is the response to the function of $f_2(x_1)$. Hence, it follows that
\begin{equation}
	f_1^{\mathrm{L}}(x_2) := \mathrm{argmax}_{x_1} u_1(x_1,f_2(x_1)).
\label{eq2-01}
\end{equation}
An obvious difference between the L- and B-responses lies in the form of the recognized player's action. Recall that in the B-response, 1's strategy is given under the image that the other's action is independent of his/her own action (see Eq.~\ref{eq1-01}). On the contrary, in the L-response, 1's strategy is given by the learning that the other's action depends on his/her own action (see Eq.~\ref{eq2-01}). Therefore, the L-response is independent of $x_2$, while the B-response depends on $x_2$.

We now consider a situation in which player 1 one-sidedly recognizes 2's intention. In this case, player 1 (2) makes the L- (B-) response. The crossing of these functions is defined as the LB equilibrium $(x_1^{\mathrm{LB}},x_2^{\mathrm{LB}})$, which is given by
\begin{equation}
\begin{split}
	& x_1^{\mathrm{LB}} = f_1^{\mathrm{L}}(x_2^{\mathrm{LB}}), \\
	& x_2^{\mathrm{LB}} = f_2^{\mathrm{B}}(x_1^{\mathrm{LB}}).
\end{split}
\label{eq2-02}
\end{equation}

Then, player $i$'s payoff is defined as $u_i^{\mathrm{LB}}:=u_i(x_1^{\mathrm{LB}},x_2^{\mathrm{LB}})$. In the same way, the BL equilibrium is defined as the crossing point between the B-response of player 1 and the L-response of player 2. In the duopoly game to be discussed later, the LB (BL) equilibrium is known as the ``Stackelberg equilibrium'' \cite{Stackelberg2011}, while in general games, it belongs to ``sub-game perfect equilibria'' \cite{Selten1965, Selten1975}. Here, we use the term the Stackelberg equilibrium in any games. Therefore, one-sided recognition means a transition from the Nash equilibrium to the Stackelberg equilibrium.

We now study some of the general properties of such one-sided recognition. First, a player does not lose any benefit by learning the other's B-response; in other words, $u_1^{\mathrm{LB}}\ge u_2^{\mathrm{BB}}$ holds. This is easily proven as
\begin{equation}
	\begin{array}{ll}
		u_1^{\mathrm{LB}} & =\mathrm{max}_{x_1}u_1(x_1,f_2^{\mathrm{B}}(x_1)) \\
		& \ge u_1(x_1^{\mathrm{BB}},f_2^{\mathrm{B}}(x_1^{\mathrm{BB}})) \\
		& =u_1^{\mathrm{BB}}. \\
	\end{array}
\label{eq2-29}
\end{equation}
This inequality is understood as follows: as the player adopting B-response chooses the strategy depending on the other, the other player can take advantage of the other's strategy and shift the equilibrium point (i.e., $x_1^{\mathrm{LB}}$ or $x_2^{\mathrm{BL}}$) one-sidedly, in order to get more payoff.  (Note that the Zero Determinant strategy by Press and Dyson \cite{Press2012} in prisoner's dilemma game, adopts a similar strategy, as the optimization strategy of one player itself is taken advantage by the other to increase the payoff.

Second, we obtain a necessary and sufficient condition for a recognizing player to increase his/her payoff. When player 1 makes the L-response, 1 refers to 2's strategy. In other words, how 2's action changes depends on 1's action. Thus, 1's action deviates from the BB equilibrium if 2's strategy function has a nonzero gradient around the BB equilibrium. Considering the case when the LB and BB equilibria are achieved within the interior of the possible range of players' actions $[x_{\mathrm{min}},x_{\mathrm{max}}]$ (i.e., at $x_{\mathrm{min}}<x<x_{\mathrm{max}}$), the condition for it is given by
\begin{equation}
\begin{split}
	&\left.\frac{\partial u_1}{\partial x_2}\right|_{\mathrm{BB}}\neq 0,\\
	&\left.\frac{\partial^2 u_2}{\partial x_1\partial x_2}\right|_{\mathrm{BB}}\neq 0.
\end{split}
\label{eq2-30}
\end{equation}
The condition for player 2's L-response is obtained in the same way. (See the Supplementary Data for the detailed calculation.)

In example 1, Eq.~\ref{eq2-30} is satisfied as long as the abilities of the two players are not equal, as discussed below. In addition, in example 2, Eq.~\ref{eq2-30} always holds (see the Supplementary Data). By contrast, in example 3, the prisoner's dilemma game, both LB and BB lie on the boundary of the range of actions $[x_{\mathrm{min}},x_{\mathrm{max}}]$ players can take. In this example, $x_i^{\mathrm{BB}} =x_i^{\mathrm{LB}}=x_i^{\mathrm{BL}}$ holds (see the Supplementary Data). Then, the learning discussed below cannot change the action nor strategy function. Hence, we no longer discuss this example.

The above result is interpreted by the relationship between both players' strategies. For the fixed 2's strategy $f_2(x_1)$, 1's strategy $f_1(x_2)$ enables him/herself to realize benefit $u_1^{\mathrm{eq}}$. Therefore, the condition that 1's strategy is optimal and is not changed by the other's strategy is given by
\begin{equation}
	u_1^{\mathrm{eq}}=\mathrm{max}_{x_1}u_1(x_1,f_2(x_1)).
\label{eq2-30_1}
\end{equation}
If the same equation for player 2 holds, the set of strategy functions is in the equilibrium. We define this as the ``function equilibrium.''

As illustrated in the following two examples, the function equilibrium is not satisfied at the BB equilibrium in general because the function of the other player imagined by one player does not agree with the real function of the other. In the B-response, the player imagines that the other's strategy function is constant, and he/she chooses his/her strategy accordingly. When both players make B-responses, however, the function of each player is no longer constant in contrast to the assumption for the B-response. Therefore, both players still gain an advantage by learning the other's strategy function.

On the contrary, there is no such disagreement at the LB or BL equilibria, where the L-response player imagines that the other's action can be changeable depending on the learning side's action, and as a result, the real strategy function is made to be constant. Thus, the real and imagined straetgy functions are consistent with each other. Then, there are no more advantages of learning the other's strategy for both players, and the function equilibrium is satisfied.

Next, we consider whether the learned side increases or decreases his/her payoff. Let us consider the ``competitive'' case in which an increase in $x_i$ leads to disutility for the other as is given by $\partial u_1/\partial x_2,\partial u_2/\partial x_1<0$. Indeed, a few nontrivial games satisfy such a relation, as discussed in these two examples. In this case, if player 1 is more competitive owing to recognition $(x_1^{\mathrm{LB}}>x_1^{\mathrm{BB}})$, the following relationship is satisfied:
\begin{equation}
	\begin{array}{ll}
		u_2^{\mathrm{LB}} &=u_2(x_1^{\mathrm{LB}},x_2^{\mathrm{LB}})\\
		&<u_2(x_1^{\mathrm{BB}},x_2^{\mathrm{LB}})\\
		&\le \mathrm{max}_{x_2}u_2(x_1^{\mathrm{BB}},x_2)\\
		&=u_2^{\mathrm{BB}}.
	\end{array}
\label{eq2-30_2}
\end{equation}
Then, the learned player is proven to receive a payoff below the BB equilibrium.

On the contrary, if player 1 is less competitive $(x_1^{\mathrm{LB}}>x_2^{\mathrm{BB}})$, we get $u_2^{\mathrm{BB}}<u_2^{\mathrm{LB}}$ in the same way. In addition, we can deal with another case, for example, a cooperative situation $\partial u_1/\partial x_2,\partial u_2/\partial x_1>0$ by inverting each player's action $x_1\rightarrow -x_1,x_2\rightarrow -x_2$. For example, the public goods game \cite{Hardin2009} belongs to this type.

Here, in contrast to earlier studies, we consider not only the equilibrium set of actions but also the functions of the players to achieve it,　based on the recognition of the other's intention. Accordingly, the function equilibrium that deviates from the Nash equilibrium is introduced. We explicitly calculate the BB, LB, and BL equilibria in specific examples.

\subsection{Example 1: resource competition game}
As an example of the BL and LB responses formulated above, we consider a ``resource competition'' game. In this game, both players $i\in\{ 1,2\}$ pay cost $x_i\ge 0$ to compete for a restricted resource with the total amount of unit one. Each player's reward, defined as the distributed resource, is proportional to the paid cost. Here, the efficiency to get resource per cost is given by $r_i$. Each player's payoff $u_i$ is defined as the difference between the reward and cost, so that
\begin{equation}
\begin{split}
	& u_1(x_1,x_2) := \frac{r_1x_1}{r_1x_1+r_2x_2} - x_1, \\
	& u_2(x_1,x_2) := \frac{r_2x_2}{r_1x_1+r_2x_2} - x_2.
\end{split}
\label{eq2-17}
\end{equation}
We assume that the abilities of the players differ, meaning that $r_1\ge r_2$. Without loss of generality, $r_2$ is set at 1, and we take $r_1\equiv r\le 1$. When $r=1$, the abilities of the two players are identical, while $r>1$ means that player 1 is superior to 2. This game is a continuous version of the hawk dove game \cite{Smith1973, Smith1988}. In addition, this continuous game was recently applied to hierarchical competition \cite{Fujimoto2017}.

% Figure 1
\begin{figure}[htb]
\begin{center}
\includegraphics[width=0.7\linewidth]{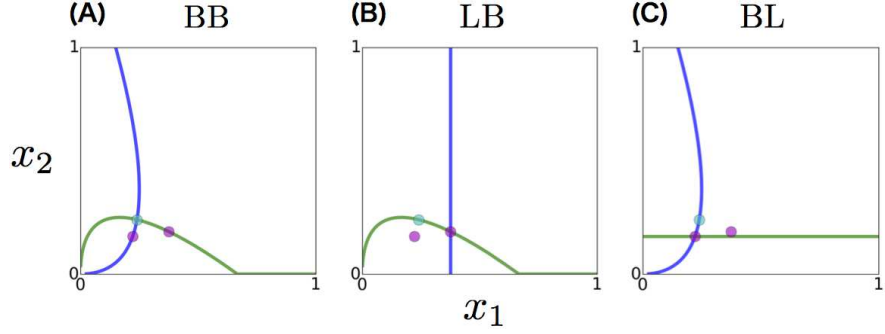}
\caption{The BB (left), LB (center), and BL (right) equilibria in the case of $r=1.5$ in a resource competition game. For all three figures, the X-axis (Y-axis) indicates player 1's (2's) action, denoted by $x_1$ ($x_2$). The blue (green) line indicates 1's (2's) intention $f_1(x_2)$ ($f_2(x_1)$) at the equilibrium. The cyan (magenta) dots indicate the Nash (Stackelberg) equilibria, respectively.}
\label{pic01}
\end{center}
\end{figure}

In this resource competition game, each player's B-response is given (see Fig.~\ref{pic01}-A) by
\begin{equation}
\begin{split}
	& f_1^{\mathrm{B}}(x_2) = \mathrm{max}(\sqrt{x_2/r}-x_2/r,0), \\
	& f_2^{\mathrm{B}}(x_1) = \mathrm{max}(\sqrt{rx_1}-rx_1,0).
\end{split}
\label{eq2-18}
\end{equation}
From Eq.~\ref{eq1-00}, we get the set of actions at the BB (Nash) equilibrium as the crossing of the strategy functions.

Next, we consider the case in which only player 1 learns 2's B-response. Player 1's L-response is given (see Fig.~\ref{pic01}-B) by
\begin{equation}
	f_1^{\mathrm{L}}(x_2) = \left\{ \begin{array}{ll}
		r/4 & (1 \le r < 2) \\
		1/r & (2 \le r) \\
	\end{array} \right.
\label{eq2-20}
\end{equation}
Then, the LB equilibrium is different from the BB equilibrium (compare Fig.~\ref{pic01}-B with A).

We now study how these two are different (see the Supplementary Data for the detailed calculation). From Fig.~\ref{pic01}, we get $x_1^{\mathrm{LB}}>x_1^{\mathrm{BB}}$, $x_2^{\mathrm{LB}}<x_2^{\mathrm{BB}}$, $u_1^{\mathrm{LB}}>u_1^{\mathrm{BB}}$, and $u_2^{\mathrm{LB}}<u_2^{\mathrm{BB}}$. These equations indicate that owing to the superior player's one-way learning, he/she increases his/her cost but increases his/her payoff, while the other player decreases his/her payoff while decreasing his/her cost.

On the contrary, when player 2 one-sidedly learns 1's B-response, 2's L-response is given (see Fig.~\ref{pic01}-C) by
\begin{equation}
	f_2^{\mathrm{L}}(x_1)=1/(4r).
\label{eq2-22}
\end{equation}
In this case, we get $x_1^{\mathrm{BL}}<x_1^{\mathrm{BB}}$, $x_2^{\mathrm{BL}}<x_2^{\mathrm{BB}}$, and $u_1^{\mathrm{BL}}>u_1^{\mathrm{BB}}$, $u_2^{\mathrm{BL}}>u_2^{\mathrm{BB}}$ as shown in Fig.~\ref{pic01} (see the Supplementary Data for the detailed calculation). Hence, both the players decrease their costs and increase their payoffs owing to the one-way learning of the inferior player in contrast to that of the superior player.

The LB and BL equilibria correspond to the classical Stackelberg ones \cite{Stackelberg2011, Selten1965, Selten1975}. In particular, LB indicates an equilibrium for a situation that 1 firstly determines his/her action and 2 follows it given the information on 1's action. BL indicates the converse situation. Note again that we focus not only on the crossing equilibrium but also on which pair of strategy functions is achieved in the equilibrium. Thus, we here call the Stackelberg equilibrium the LB or BL equilibrium in the same way that we call the Nash equilibrium the BB equilibrium.

The superior player's one-way learning results in exploitation by gaining more benefit by increasing its own cost, while the inferior player's learning results in cooperation by decreasing its own cost. This is interpreted as follows. First, the cost a player should pay depends on the other's cost. A player would not need to pay so much when the other's cost is too small because the player would monopolize most resources by paying not so much cost. On the contrary, if the other's cost is too large, the player would not pay much cost either because one should pay too much cost to obtain more resources. Therefore, a player's optimal cost is maximal when the other pays a moderate cost (see Fig.~\ref{pic01}-A). Second, in the BB equilibrium, both players pay a moderate cost. Hence, no matter whether the learning player is superior or not, the player has to repress the other's cost to gain more benefits. How to repress the other's cost, however, depends on whether the player is superior or not. The superior player increases his/her cost and forces the inferior one to give up competition. Therefore, the former exploits the latter by learning the other's strategy. On the contrary, to gain a higher payoff, the inferior player decreases his/her cost and relaxes the competition. Therefore, the learning player cooperates with the learned one.

How the learning and learned players' payoffs change depends on the type of game. Below, we discuss an alternative example, namely a duopoly game, in which competition always persists to the point that an increase in the payoff of the learning player always decreases the other's payoff.

\subsection{Example 2: duopoly game}
In a duopoly game, two companies $i\in\{1,2\}$, which separately supply products, compete for a limited market. The more products supplied in this limited market, the cheaper their prices are. Here, player $i$'s action $x_i$ is the number of products he/she supplies. We assume that the price is $\max(0,p-x_1-x_2)$, where $p$ represents the maximal price. In addition, player $i$'s cost of supplying products is assumed to be $c_i$. Thus, player $i$'s payoff is given by
\begin{equation}
\begin{split}
	u_1(x_1,x_2):=x_1(\max(0,p-x_1-x_2)-c_1),\\
	u_2(x_1,x_2):=x_2(\max(0,p-x_1-x_2)-c_2).
\end{split}
\label{eq2-23}
\end{equation}
Here, we assume $c_1 \le c_2$ without loss of generality. In other words, player 1 is superior to 2.

\begin{figure}[htb]
\begin{center}
\includegraphics[width=0.7\linewidth]{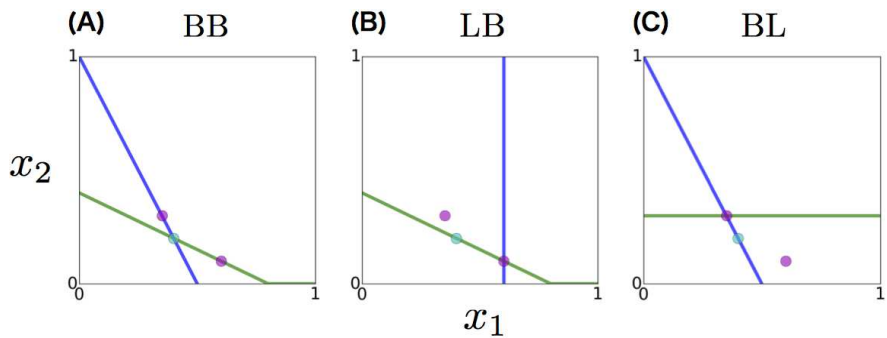}
\caption{The BB (left), LB (center), and BL (right) equilibria in the case of $(c_1,c_2)=(0,0.2)$ in a resource competition game. For all three figures, the X-axis (Y-axis) indicates player 1's (2's) action, denoted by $x_1$ ($x_2$). The blue (green) line indicates 1's (2's) strategy function $f_1(x_2)$ ($f_2(x_1)$) at the equilibrium. The cyan (magenta) dots indicate the Nash (Stackelberg) equilibria.}
\label{pic02}
\end{center}
\end{figure}

The strategy functions in the BB, LB, and BL equilibria are plotted in Figs.~\ref{pic02}-A, B, and C, respectively. Accordingly, from the crossing points, we obtain $u_1^{\mathrm{LB}}>u_1^{\mathrm{BB}}$ and $u_2^{\mathrm{LB}}<u_2^{\mathrm{BB}}$. These equations indicate that the superior company's learning results in the exploitation of the learned one, as in the resource competition game. On the contrary, we get $u_1^{\mathrm{BL}}<u_1^{\mathrm{BB}}$ and $u_2^{\mathrm{BL}}>u_2^{\mathrm{BB}}$. In contrast to the resource competition game, the inferior company's learning also results in exploitation. (See the Supplementary Data for the detailed calculation.)

% Figure 2

Although both the resource competition and the duopoly games are categorized as competitive ($\partial u_1/\partial x_2, \partial u_2/\partial x_1<0$), the change in the learned player's payoff differs between the two. As has already been explained, this difference depends on whether the learning side is more or less competitive according to one-way learning (see Eq.~\ref{eq2-30_2}). Each player's motivation to change his/her competitiveness is now discussed based on the following dynamic process of learning.

\section{Functional dynamics of strategies}
So far, we first considered the BB equilibrium in which both players have no recognition of the other's intention. Then, we introduced the LB and BL equilibria in which a player one-sidedly recognizes the intention of the other. When one recognizes the other but not vice versa, the recognizing player has an advantage. However, such one-way recognition rarely appears because both of the players usually try to recognize each other's intention. Another problem in the L-response is that one player knows the other's intention perfectly in a one-shot game, while players usually shape the image of the other successively through the iteration of games. If both players recognize the other's intention gradually, the LB or BL (Stackelberg) equilibrium is no longer achieved. Instead, a set of actions that are not discussed by previous studies can be achieved, as shown below.

To represent such a gradual recognizing process, we assume that each player $i$ learns the other's strategy at a rate of $\epsilon_i$. In this case, each player's strategy function $f_1(x_2)$ and $f_2(x_1)$ changes depending on the other's one as
\begin{equation}
\begin{split}
	f_1(x_2)[t+1]=\mathrm{argmax}_{x_1}u_1(x_1,\epsilon_1f_2(x_1)[t]+(1-\epsilon_1)x_2), \\
	f_2(x_1)[t+1]=\mathrm{argmax}_{x_2}u_2(\epsilon_2f_1(x_2)[t]+(1-\epsilon_2)x_1,x_2).
\end{split}
\label{eq3-01}
\end{equation}

For $\epsilon_1=0$, $f_1$ in the one-shot game corresponds to the B-response (see Eq.~\ref{eq1-01}); for $\epsilon_1=1$, $f_1$ corresponds to the L-response (see Eq.~\ref{eq2-01}). In addition, when at least one player makes the B-response, both players' strategies in the equilibrium are given as fixed functions, as already mentioned. In the present case with $\epsilon_1,\epsilon_2>0$, however, it is necessary to consider the functional dynamics, where both players change their strategy functions by learning the other's strategy function. Therefore, we add the time variable $t$.

Eq.~\ref{eq3-01} represents the functional dynamics \cite{Kataoka2000, Kataoka2001, Kataoka2003}, where the change in time depends on the function rather than the dynamic systems of state variables of a finite dimension (for example, in dynamical-systems game \cite{Sato2002}). Hence, we need to solve the dynamics of infinite dimensions.

We now analyze the equilibrium state of Eq.~\ref{eq3-01}. In the following, we assume that there exist a pair of fixed-point functions as an equilibrium state of the functional dynamics, which is denoted by $f_1^{*}(x_2), f_2^{*}(x_1)$ satisfying
\begin{equation}
\begin{split}
	f_1^{*}(x_2) = \mathrm{argmax}_{x_1} u_1(x_1,\epsilon_1f_2^{*}(x_1)+(1-\epsilon_1)x_2), \\
	f_2^{*}(x_1) = \mathrm{argmax}_{x_2} u_2(\epsilon_2f_1^{*}(x_2)+(1-\epsilon_2)x_1,x_2).
\end{split}
\label{eq3-01_3}
\end{equation}
As demonstrated later numerically, fixed-point functions are reached in various games. To study the behavior near the equilibrium, we derive a crossing point of the equilibrium functions and its neighborhood. By assuming the continuity of the functions around the crossing point, we expand the equilibrium functions as
\begin{equation}
\begin{split}
	f_1^{*}(x_2) = x_1^{\mathrm{eq}*} + a_1^{*}(x_2-x_2^{\mathrm{eq}*}), \\
	f_2^{*}(x_1) = x_2^{\mathrm{eq}*} + a_2^{*}(x_1-x_1^{\mathrm{eq}*}).
\end{split}
\label{eq3-01_2}
\end{equation}
By substituting Eq.~\ref{eq3-01_2} into Eq.~\ref{eq3-01_3}, we get the first-order term as
\begin{equation}
\begin{split}
	\left.\frac{\partial}{\partial x_1}\left(\left.u_1(x_1,x_2)\right|_{x_2=\hat{x_2}}\right)\right|_{\mathrm{eq}*}=0, \\
	\left.\frac{\partial}{\partial x_2}\left(\left.u_2(x_1,x_2)\right|_{x_1=\hat{x_1}}\right)\right|_{\mathrm{eq}*}=0.
\end{split}
\label{eq3-01_4}
\end{equation}
Then, we also get the second-order term as
\begin{equation}
\begin{split}
	\left.\frac{\partial}{\partial x_2}\left(\left.\frac{\partial}{\partial x_1}\left(\left.u_1(x_1,x_2)\right|_{x_2=\hat{x_2}}\right)\right|_{x_1=f_1^{*}(x_2)}\right)\right|_{\mathrm{eq}*}=0, \\
	\left.\frac{\partial}{\partial x_1}\left(\left.\frac{\partial}{\partial x_2}\left(\left.u_2(x_1,x_2)\right|_{x_1=\hat{x_1}}\right)\right|_{x_2=f_2^{*}(x_1)}\right)\right|_{\mathrm{eq}*}=0,
\end{split}
\label{eq3-01_5}
\end{equation}
where $\hat{x_1}=:\epsilon_2f_1^{*}(x_2)+(1-\epsilon_2)x_1$ and $\hat{x_2}:=\epsilon_1f_2^{*}(x_1)+(1-\epsilon_1)x_2$. Here, Eq.~\ref{eq3-01_4} indicates that the crossing point $(x_1^{\mathrm{eq}*},x_2^{\mathrm{eq}*})$ satisfies the optimization condition for the other's function, while Eq.~\ref{eq3-01_5} is a consequence of the {\bf fixed-point functions}, indicating that the set of equilibrium functions $(f_1^{*}(x_2),f_2^{*}(x_1))$ satisfies the optimization condition in the neighborhood of the crossing point. To compute the player's equilibrium payoff given by $(x_1^{\mathrm{eq}*},x_2^{\mathrm{eq}*})$, the above calculation of $f_1^{*}$ and $f_2^{*}$ is thus sufficient.

\begin{figure}[htb]
\begin{center}
\includegraphics[width=0.7\linewidth]{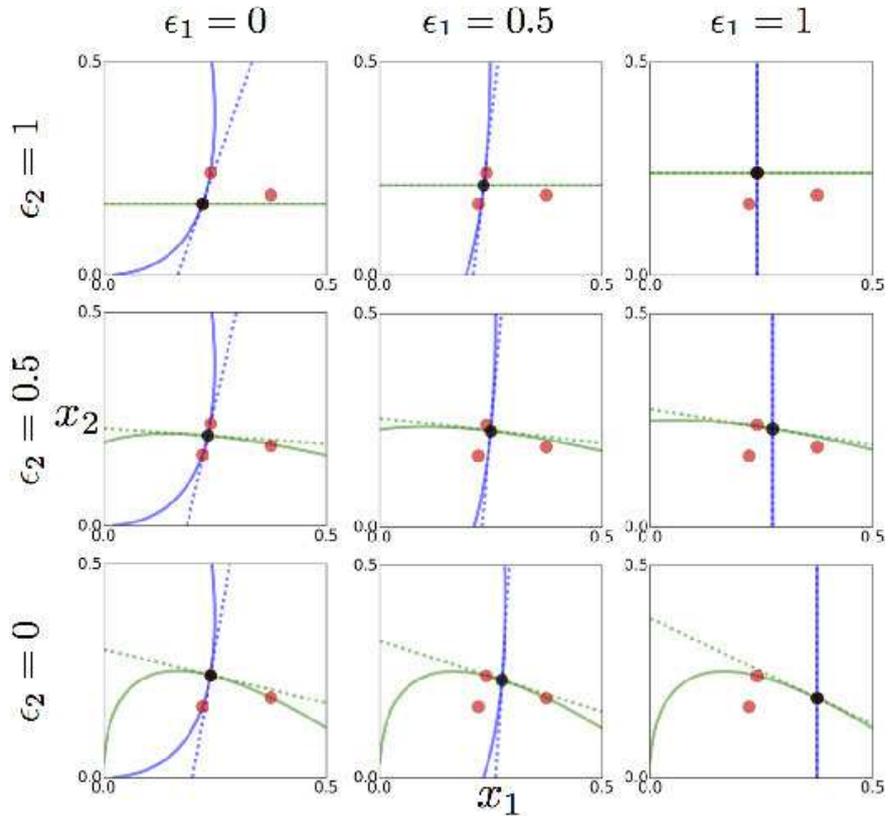}
\caption{Equilibrium state of mutual learning in a resource competition game. For all nine figures, the X-axis (Y-axis) indicates player 1's (2's) action, denoted by $x_1$ ($x_2$). The blue (green) line indicates 1's (2's) strategy function and the solid (broken) line indicates the simulated (analytical) solution. The yellow circle is the crossing point and the red circles are the Nash and Stackelberg equilibria to allow for a comparison of positions. The left, center, and right figures are respectively the cases of $\epsilon_1=0,0.5,1$, while the upper, center, and lower figures are respectively the cases of $\epsilon_2=0,0.5,1$.}
\label{pic03}
\end{center}
\end{figure}

As an extreme case, we consider $\epsilon_i=1$, in which both players perfectly recognize the other's intention. The fixed point in this case is the LL equilibrium. Here, both players make L-responses, which are constant for the other's action. Therefore, from Eq.~\ref{eq2-01}, the achieved actions at the LL equilibrium correspond to those at BB (i.e., the Nash equilibrium), namely $x_i^{\mathrm{LL}}=x_i^{\mathrm{BB}}$ (see the upper right of Fig.~\ref{pic03}). Note that the equilibrium points are identical, whereas they are different in the functional dynamics. Indeed, from Eqs.~\ref{eq3-01_4} and \ref{eq3-01_5}, we can confirm that $x_i^{\mathrm{eq}*}$ is equal for LL and BB, while $a_i^{*}$ is not (see the Supplementary Data for the detailed calculation). Owing to this inequality in $a_i^{*}$, the function equilibrium holds in LL, but not in BB: In the functional dynamics with $\epsilon_i=0$, both players' strategy functions are constant in LL, while for $\epsilon_i=1$, they are not.

% Figure 3

It may be disappointing that the classical Nash equilibrium is achieved as the LL equilibrium in the game between the learning players. As explained later, however, this LL equilibrium is rarely achieved; indeed, in many cases, novel equilibria are achieved according to the functional dynamics.

Below, we discuss some specific examples for the functional dynamics.

\subsection{Example 1: resource competition game}
\begin{figure}[htb]
\begin{center}
\includegraphics[width=0.7\linewidth]{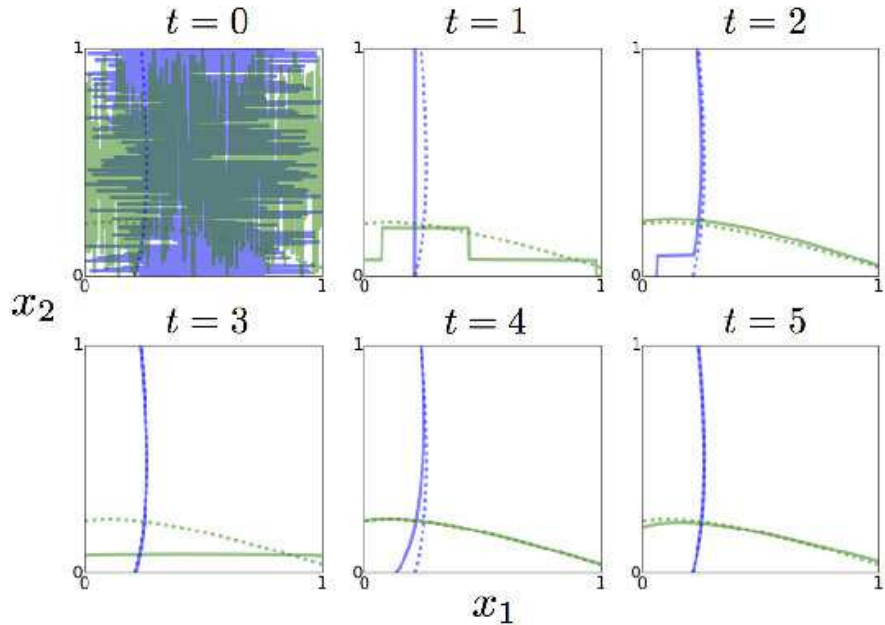}
\caption{An example of the functional dynamics in a resource competition game for $r=1.5$, $(\epsilon_1,\epsilon_2)=(0.5,0.5)$. For all figures, the X-axis (Y-axis) indicates player 1's (2's) action, denoted by $x_1$ ($x_2$). The solid blue (green) line indicates 1's (2's) strategy function until $t=6$. The broken lines indicate both players' ones in the equilibrium, achieved at $t=20$.}
\label{pic04}
\end{center}
\end{figure}

We again consider the resource competition game. From Eqs.~\ref{eq3-01_4} and \ref{eq3-01_5}, we get the set of equilibrium actions $(x_1^{\mathrm{eq}*}, x_2^{\mathrm{eq}*})$ and the set of equilibrium gradients around them $(a_1^{*}, a_2^{*})$ as given by
\begin{equation}
\begin{split}
	& x_1^{\mathrm{eq}*} = \frac{1}{r+\epsilon_1a_2^{*}} \left\{ \sqrt{r(x_2^{\mathrm{eq}*}-\epsilon_1a_2^{*}x_1^{\mathrm{eq}*})}-(x_2^{\mathrm{eq}*}-\epsilon_1a_2^{*}x_1^{\mathrm{eq}*}) \right\}, \\
	& x_2^{\mathrm{eq}*} = \frac{1}{1+r\epsilon_2a_1^{*}} \left\{ \sqrt{r(x_1^{\mathrm{eq}*}-\epsilon_2a_1^{*}x_2^{\mathrm{eq}*})}-r(x_1^{\mathrm{eq}*}-\epsilon_2a_1^{*}x_2^{\mathrm{eq}*}) \right\}, \\
	& a_1^{*} = \frac{1}{r+\epsilon_1a_2^{*}} \left\{ \frac{r(1-\epsilon_1)}{2\sqrt{r(x_2^{\mathrm{eq}*}-\epsilon_1a_2^{*}x_1^{\mathrm{eq}*})}} - (1-\epsilon_1) \right\}, \\
	& a_2^{*} = \frac{1}{1+r\epsilon_2a_1^{*}} \left\{ \frac{r(1-\epsilon_2)}{2\sqrt{r(x_1^{\mathrm{eq}*}-\epsilon_2a_1^{*}x_2^{\mathrm{eq}*})}} - r(1-\epsilon_2) \right\}.
\end{split}
\label{eq3-09}
\end{equation}
We now simulate Eq.~\ref{eq3-01} and compare the simulation results with the calculation, confirming that both players' strategy functions immediately converge to fixed ones (see Fig.~\ref{pic04}). The crossing points $(x_1^{\mathrm{eq}*}, x_2^{\mathrm{eq}*})$ of these functions agree well with the above analytic estimation and the converged strategy functions in the neighborhood of the crossing points are well estimated by Eq.~\ref{eq3-01_2} with the above values $a_1^{*}$ and $a_2^{*}$. In addition, the action compared with the other's action is calculated from fixed-point function $f$, as shown in Fig.~\ref{pic03}. This indicates that the more (less) each player learns the other's strategy, the less (more) dependent the strategy function is on the other's action.

% Figure 4

\subsection{Example 2: duopoly game}
As in the resource competition game, both players' strategy functions converge to fixed functions for any $\epsilon_1$ and $\epsilon_2$. Here, recall that the learning side exploits the learned side regardless of whether the former is superior or not. This result can be applied to the case with continuous learning degrees $\epsilon_1$ and $\epsilon_2$. The larger $\epsilon_1$, the larger (smaller) 1's (2's) payoff is, with larger exploitation (see the Supplementary Data for the details).

\section{Dynamics of the degree of learning}
Thus far, the learning degree $\epsilon_i$ has been given and fixed. Thus, for each player, the case with $\epsilon_i=1$ would be the better one for receiving a higher payoff. Each player, however, can change the degree to which he/she learns the other's strategy. Initially, each player may not care about the other, and he/she learns the other's strategy more through the repeated game. In the following, we consider this temporal evolution in the degree of learning, $\epsilon_1,\epsilon_2$. Here, assuming that the other's strategy function is fixed, each player tries to increase his/her payoff by changing his/her learning degree. Therefore, the dynamics of both players' learning degrees are given by
\begin{equation}
\begin{split}
	\dot{\epsilon_1}=S_1\frac{\partial\left.u_1(x_1,\epsilon_1f_2(x_1)+(1-\epsilon_1)x_2)\right|_{\mathrm{eq}*}}{\partial \epsilon_1}, \\
	\dot{\epsilon_2}=S_2\frac{\partial\left.u_2(\epsilon_2f_1(x_2)+(1-\epsilon_2)x_1,x_2)\right|_{\mathrm{eq}*}}{\partial\epsilon_2},
\end{split}
\label{eq3-02}
\end{equation}
where $S_1,S_2$ is the speed with which each player optimizes the intensity of recognition. In the following, we simulate these dynamics for the introduced example and examine what equilibrium is reached.

\subsection{Example 1: resource competition game}
\begin{figure}[htb]
\begin{center}
\includegraphics[width=0.5\linewidth]{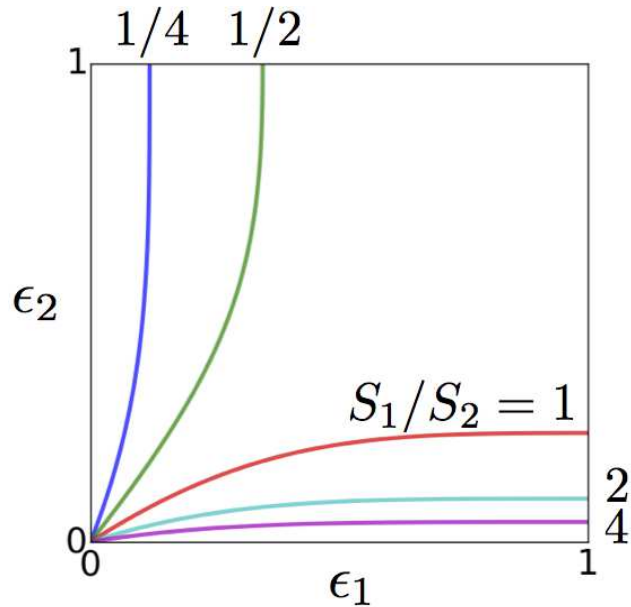}
\caption{Trajectories of a set of learning degrees $(\epsilon_1, \epsilon_2)$ for different ratios of learning speeds $S_1/S_2$ in a resource competition game. The blue, green, red, cyan, and magenta lines represent the cases of $S_1/S_2=0.25, 0.5, 1, 2, 4$, respectively. For all trajectories, the dynamics of Eq.~\ref{eq3-02} start from $(\epsilon_1,\epsilon_2)=(0,0)$.}
\label{pic05}
\end{center}
\end{figure}

\begin{figure}[htb]
\begin{center}
\includegraphics[width=0.6\linewidth]{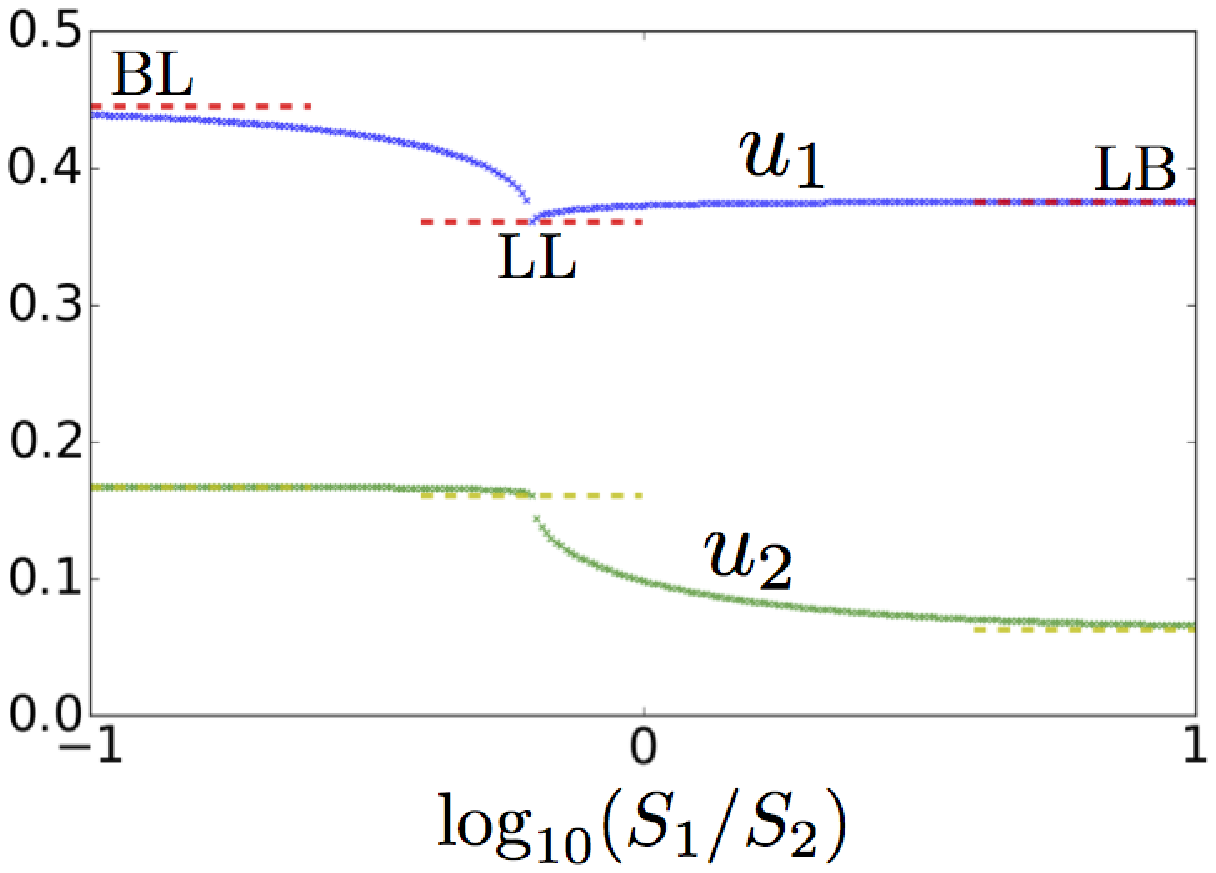}
\caption{Player 1's payoff (blue dots) and player 2's payoff (green dots) in the equilibrium state of Eq.~\ref{eq3-02} for the diverse sets of learning speeds $(S_1,S_2)$ in a resource competition game. The red (yellow) broken lines indicate player 1's (2's) payoff. The right, middle, and left broken lines indicate the payoffs at the LB, LL, and BL equilibria, respectively.}
\label{pic06}
\end{center}
\end{figure}

Fig.~\ref{pic05} shows the dynamics of $\epsilon_1$ and $\epsilon_2$ for various sets of learning speeds $(S_1,S_2)$, while this temporal evolution in the payoff of each player according to Eq.~\ref{eq3-02} is shown in Fig.~\ref{pic06}.

% Figure 5

% Figure 6

First, the initial BB (Nash) equilibrium is unstable compared with the learning dynamics. In other words, both players are motivated to learn the other's strategy and to change their strategies accordingly because the payoff has a nonzero gradient at $(x_1,x_2)$ around the BB equilibrium point.

Second, as one player's learning is superior, the other's learning is repressed. During the evolution of learning, one player learns the other's gradient of the strategy function, and his/her strategy approaches a constant function (L-response). Each player tries to gain more benefit by repressing the other's cost. Therefore, the superior player 1 increases his/her cost as a result of his/her own learning, while the inferior player 2 decreases his/her cost.

Third, the intermediate state between the LB/BL and LL equilibria is finally achieved depending on the evolutionary speed relationship. Not only the edges (LB, BL, and LL equilibria), but also the intermediate states of LB (BL) and LL ($\epsilon_1,\epsilon_2$ is between 0 and 1) are achieved, where every player succeeds in optimizing his/her strategy function for the other's one. In other words, the function equilibrium is achieved as a result of these dynamics.

\subsection{Example 2: duopoly game}
In this case, the players also finally reach the intermediate state between the LB/BL and LL equilibria as a result of mutual learning, where the function equilibrium is satisfied. In addition, the more 1's (2's) learning speed increases, the more 1 (2) exploits the benefit from the other in the final state (see the Supplementary Data for more details.)

\section{Summary and discussion}
In this study, we introduce a new formulation for the mutual recognition of intention, which is represented as functional dynamics. In the formulation, every player can read the other's strategy function $f(x)$, which determines the action to be chosen for the other's action $x$.

As a result, we proved that both players can increase their payoffs according to their learning. The more a player learns the other's strategy, the less his/her action depends on the other's action (i.e., the function approaches a constant function). Since such a constant function does not provide any motivation to learn, the process of mutual learning stops when one player perfectly learns the other's strategy function.

The resultant function equilibrium includes and extends two kinds of well-known equilibria (Nash and Stackelberg), which are achieved in contrasting situations. In the extreme case in which one player one-sidedly learns the other, the Stackelberg equilibrium is achieved. On the contrary, when the evolutionary speeds of the learning of both players are of comparative order, the Nash equilibrium is achieved again. In contrast to previous studies, our study discusses such equilibria as the crossing points of the strategy functions. Furthermore, depending on the evolutionary speeds of learning, the intermediate equilibria between the Nash and Stackelberg equilibria are achieved. This finding may lead to understanding how the leader--follower relationship is formed in game theory according to intention recognition.

Furthermore, we also confirm that each player's payoff changes according to intention recognition in the resource competition game. The faster the superior player learns the other's strategy function, the more that player exploits the other, where the learning (learned) player's payoff increases (decreases). On the contrary, the faster the inferior person learns, the more the players cooperate, meaning that both players' payoffs increase. In the learning process, one player's positive (negative) gradient of his/her strategy function leads to a decrease (increase) in the other's competitiveness.

In this study, we only consider the two-player case. In the case of $n(\ge 3)$ players, each player's strategy function is $(n-1)$-dimensional, and the function equilibrium is more complicated. Furthermore, some equilibria regarded as neither Nash nor Stackelberg equilibria emerge. For example, we can consider the cases of one-sided learning from 1 to 2, from 2 to 3, and from 3 to 1. Such a learning loop does not appear as an equilibrium in the extensive form game \cite{Kuhn1950}. This will be discussed in future work.

In our formulation, we assume that each player's learning speed is independent of its accuracy. This assumption results in a monotonic advantage for the increase in learning speed, at least in the duopoly game. In reality, however, there is a trade-off between the accuracy of reading and speed of evolution, which provides another disadvantage for the fast evolution owing to incomplete information. Indeed, some previous studies show that such incomplete information on the other's action leads to disutility \cite{Wu1995}.

%In this study, the deviation from the Nash equilibrium does not always happen (e.g., in the prisoner's dilemma game). A previous study, however, suggests that cooperation can be achieved by deviating from the Nash equilibrium, even in the prisoner's dilemma game, because we assume both players' initial states to be the best responses. Hence, a discussion of the response function dynamics from the general initial states will be an important extension in future research.

\section{acknowledgement}
The authors would like to thank E. Akiyama, T. S. Hatakeyama and N. Takeuchi for useful discussions. This research was partially supported by Hitachi The University of Tokyo Laboratory.

% Bibliograph

%%%%%
%Supplementary Material%

\newpage

\setcounter{figure}{0}
\renewcommand{\figurename}{FIG. S}

\begin{center}
{\LARGE {\bf Supplementary Material}}
\end{center}

\section{General games}
\subsection{Condition for disagreement between the Nash and Stackelberg equilibria}
Here, we consider a condition under which a player increases his/her payoff owing to the one-sided recognition of the other's response. In the following, we assume a game satisfying the following two points: (1) the space of each player's possible action is bounded and (2) each player's payoff $u_i(x_1,x_2)$ is a twice differentiable function. Then, we consider two situations for agreement between the LB and BB equilibria: both of them exist (i) in the inside of the space of all players' actions or (ii) on the border. In the following, we discuss whether LB and BB coincide for each of these situations.

First, when the Nash and Stackelberg equilibria exist within the space of all players' actions (examples 1 and 2 satisfy this condition in reality), a set of actions at the Nash equilibrium $(x_1^{\mathrm{BB}},x_2^{\mathrm{BB}})$ satisfies
\begin{eqnarray}
	&\displaystyle\left.\frac{\partial u_1}{\partial x_1}\right|_{\mathrm{BB}}=0,
\label{eq1-01_1}\\
	&\displaystyle\left.\frac{\partial u_2}{\partial x_2}\right|_{\mathrm{BB}}=0.
\label{eq1-01_2}
\end{eqnarray}
Next, we consider a condition satisfied by the Stackelberg (LB) equilibrium. When player 1's equilibrium action $x_1^{\mathrm{LB}}$ is optimal given the recognition of the other's intention $f_2^{\mathrm{B}}$, we get
\begin{equation}
\begin{split}
	&\left.\frac{\partial u_1(x_1,f_2^{\mathrm{B}}(x_1))}{\partial x_1}\right|_{\mathrm{LB}}=0,\\
	&\Leftrightarrow\left.\frac{\partial u_1}{\partial x_1}\right|_{\mathrm{LB}}+\left.\frac{\partial f_2^{\mathrm{B}}}{\partial x_1}\right|_{\mathrm{LB}}\left.\frac{\partial u_1}{\partial x_2}\right|_{\mathrm{LB}}=0.
\end{split}
\label{eq1-02}
\end{equation}
Here, $\partial f_2^{\mathrm{B}}/\partial x_1$ is given by
\begin{equation}
\begin{split}
	&\left.\frac{\partial u_2}{\partial x_2}\right|_{x_2=f_2^{\mathrm{B}}(x_1)}=0,\\
	&\Leftrightarrow\frac{\partial^2u_2}{\partial x_1\partial x_2}+\frac{\partial f_2^{\mathrm{B}}}{\partial x_1}\frac{\partial^2u_2}{\partial x_2^2}=0.
\end{split}
\label{eq1-03}
\end{equation}
By substituting Eq.~\ref{eq1-03} into Eq.~\ref{eq1-02}, we obtain a condition for the LB equilibrium as
\begin{equation}
	\left.\frac{\partial u_1}{\partial x_1}\right|_{\mathrm{LB}}\left.\frac{\partial^2 u_2}{\partial x_2^2}\right|_{\mathrm{LB}}+\left.\frac{\partial u_1}{\partial x_2}\right|_{\mathrm{LB}}\left.\frac{\partial^2u_2}{\partial x_1\partial x_2}\right|_{\mathrm{LB}}=0.
\label{eq1-04_1}
\end{equation}
In addition, since player 2 makes a B-response, another condition for the LB equilibrium is given by
\begin{equation}
	\left.\frac{\partial u_2}{\partial x_2}\right|_{\mathrm{LB}}=0.
\label{eq1-04_2}
\end{equation}
Let us now compare the above condition (Eqs.~\ref{eq1-04_1} and \ref{eq1-04_2}) for the LB equilibrium with that of the BB equilibrium (Eqs.~\ref{eq1-01_1}, \ref{eq1-01_2}). First, the second condition is common between the two. Next, the first condition is different as long as the second term in Eq.~\ref{eq1-04_1} is nonzero. Hence, the condition for the mismatch between LB and BB is given by
\begin{equation}
\begin{split}
	&\left.\frac{\partial u_1}{\partial x_2}\right|_{\mathrm{BB}}\neq 0,\\
	&\left.\frac{\partial^2 u_2}{\partial x_1\partial x_2}\right|_{\mathrm{BB}}\neq 0.
\end{split}
\label{eq1-05}
\end{equation}
The condition for disagreement between BL and BB is given in the same way. In examples 1 and 2, we discuss this condition concretely.

\begin{figure}[htb]
\begin{center}
\includegraphics[width=0.4\linewidth]{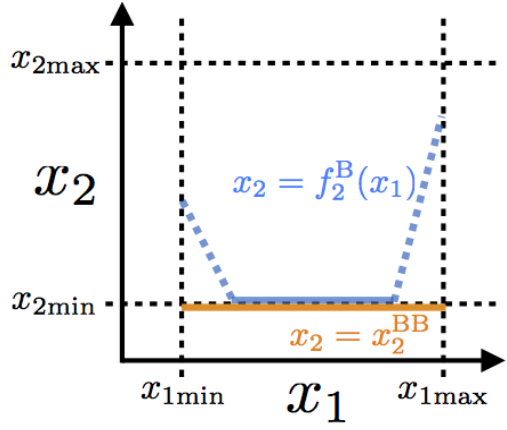}
\caption{Player 2's response functions when LB and BB exist on $x_2=x_2^{\mathrm{min}}$. The blue (green) line indicates player 2's response in LB $x_2=f_2^{\mathrm{B}}(x_1)$ (BB $x_2=x_2^{\mathrm{BB}}$). From the above assumption, the LB equilibrium can only exist on the blue solid line. Thus, the existence region of the LB equilibrium (the first line in Eq.~\ref{eq1-05_5}) is included by that of the BB equilibrium (the second line in Eq.~\ref{eq1-05_5}).}
\label{picS05}
\end{center}
\end{figure}

Second, we also consider a case in which the Nash and Stackelberg equilibria exist on the border of the space of all players' actions. In particular, we consider a situation that the LB and BB equilibria are on the border of player 2's action, in other words, $x_2^{\mathrm{LB}}=x_2^{\mathrm{BB}}=\mathrm{inf}x_2$ or $\mathrm{sup}x_2$ holds. In this case, the recognizer's payoff satisfies
\begin{equation}
\begin{split}
	u_1^{\mathrm{LB}}&=\mathrm{max}_{x_1}u_1(x_1,f_2^{\mathrm{B}}(x_1))\\
	&\le\mathrm{max}_{x_1}u_1(x_1,x_2^{\mathrm{BB}})\\
	&=u_1^{\mathrm{BB}}.
\end{split}
\label{eq1-05_5}
\end{equation}
Here, from the assumption that the LB equilibrium exists on the border of player 2's action, we can derive the second line from the first one (see Fig.~S\ref{picS05}). In addition, as shown in the main manuscript, player 1 obtains no advantage by recognizing the other's B-response; in other words, $u_1^{\mathrm{LB}}\ge u_1^{\mathrm{BB}}$ holds. Thus, we get $u_1^{\mathrm{LB}}=u_1^{\mathrm{BB}}$ and then prove $x_1^{\mathrm{LB}}=x_1^{\mathrm{BB}}$. Putting this in other terms, $x_1^{\mathrm{LB}}\neq x_1^{\mathrm{BB}}$ and $x_2^{\mathrm{LB}}=x_2^{\mathrm{BB}}=\mathrm{inf}x_2,\mathrm{sup}x_2$ are incompatible. In example 3, we see a concrete illustration of the overlap between the LB and BB equilibria.

% Figure S05

\subsection{Solution of the functional dynamics}
Here, we derive the following functional dynamics:
\begin{equation}
\begin{split}
	&f_1^{*}(x_2) = \mathrm{argmax}_{x_1} u_1(x_1,\epsilon_1f_2^{*}(x_1)+(1-\epsilon_1)x_2) \\
	&f_2^{*}(x_1) = \mathrm{argmax}_{x_2} u_2(\epsilon_2f_1^{*}(x_2)+(1-\epsilon_2)x_1,x_2)
\end{split}
\label{eq1-06}
\end{equation}
First, we expand the equilibrium functions around the crossing points as
\begin{equation}
\begin{split}
	&f_1^{*}(x_2) = a_1^{*}(x_2-x_2^{\mathrm{eq}*}) + x_1^{\mathrm{eq}*}, \\
	&f_2^{*}(x_1) = a_2^{*}(x_1-x_1^{\mathrm{eq}*}) + x_2^{\mathrm{eq}*}.
\end{split}
\label{eq1-07}
\end{equation}

In the following, we calculate the solution of player 1's function $f_1^{*}(x_2)$, which is given by
\begin{equation}
	\left.\frac{\partial}{\partial x_1}\left(u_1(x_1,\epsilon_1f_2^{*}(x_1)+(1-\epsilon_1)x_2)\right)\right|_{x_1=f_1^{*}(x_2)}=0
\label{eq1-08}
\end{equation}
in $x_2 \simeq x_2^{\mathrm{eq}*}$. Thus, we obtain the first-order term as
\begin{equation}
\begin{split}
	&\left.\left.\frac{\partial}{\partial x_1}\left(u_1(x_1,\epsilon_1f_2^{*}(x_1)+(1-\epsilon_1)x_2)\right)\right|_{x_1=f_1^{*}(x_2)}\right|_{x_2=x_2^{\mathrm{eq}*}}=0\\
	&\Leftrightarrow\left.\frac{\partial}{\partial x_1}\left(u_1(x_1,\epsilon_1f_2^{*}(x_1)+(1-\epsilon_1)x_2)\right)\right|_{\mathrm{eq}*}\\
	&\Leftrightarrow\left.\frac{\partial u_1}{\partial x_1}\right|_{\mathrm{eq}*}+\epsilon_1a_2^{*}\left.\frac{\partial u_1}{\partial x_2}\right|_{\mathrm{eq}*}=0.
\end{split}
\label{eq1-09}
\end{equation}
Here, $\left.\right|_{\mathrm{eq}*}$ is an operation of $(x_1,x_2)=(x_1^{\mathrm{eq}*},x_2^{\mathrm{eq}*})$. Then, we also obtain the second-order term as
\begin{equation}
\begin{split}
	&\frac{\partial}{\partial x_2}\left.\left(\left.\frac{\partial}{\partial x_1}\left(u_1(x_1,\epsilon_1f_2^{*}(x_1)+(1-\epsilon_1)x_2)\right)\right|_{x_1=f_1^{*}(x_2)}\right)\right|_{x_2=x_2^{\mathrm{eq}*}}=0\\
	&\Leftrightarrow a_1^{*}\left(\left.\frac{\partial^2u_1}{\partial x_1^2}\right|_{\mathrm{eq}*}+2\epsilon_1a_2^{*}\left.\frac{\partial^2u_1}{\partial x_1\partial x_2}\right|_{\mathrm{eq}*}+(\epsilon_1a_2^{*})^2\left.\frac{\partial^2u_1}{\partial x_2^2}\right|_{\mathrm{eq}*}\right)+(1-\epsilon_1)\left(\left.\frac{\partial^2u_1}{\partial x_1\partial x_2}\right|_{\mathrm{eq}*}+\epsilon_1a_2^{*}\left.\frac{\partial^2u_1}{\partial x_2^2}\right|_{\mathrm{eq}*}\right)=0.
\end{split}
\label{eq1-10}
\end{equation}

\section{Example 1: resource competition game}
In this game, both players' payoffs are defined by
\begin{equation}
\begin{split}
	& u_1(x_1,x_2) := \frac{rx_1}{rx_1+x_2} - x_1, \\
	& u_2(x_1,x_2) := \frac{x_2}{rx_1+x_2} - x_2.
\end{split}
\label{eq2-01}
\end{equation}

\subsection{BB equilibrium}
Each player's B-response is calculated as follows:
\begin{equation}
\begin{split}
	& \begin{array}{ll}
		f_1^{\mathrm{B}}(x_2) & :=\mathrm{argmax}_{x_1\ge 0}u_1(x_1,x_2) \\
		& =\max(\sqrt{x_2/r}-(x_2/r),0) \\
	\end{array},\\
	& \begin{array}{ll}
		f_2^{\mathrm{B}}(x_1) & :=\mathrm{argmax}_{x_2\ge 0}u_2(x_1,x_2) \\
		& =\max(\sqrt{x_1r}-(x_1r),0) \\.
	\end{array}
\end{split}
\label{eq2-02}
\end{equation}

Then, the crossing actions at the BB equilibrium are given as
\begin{equation}
\begin{split}
	&\left\{\begin{array}{l}
		x_1^{\mathrm{BB}}=\max(\sqrt{x_2^{\mathrm{BB}}/r}-(x_2^{\mathrm{BB}}/r),0)\\
		x_2^{\mathrm{BB}}=\max(\sqrt{x_1^{\mathrm{BB}}r}-(x_1^{\mathrm{BB}}r),0)\\
	\end{array}\right.\\
	&\Leftrightarrow x_1^{\mathrm{BB}}=x_2^{\mathrm{BB}}=r/(1+r)^2.
\end{split}
\label{eq2-03}
\end{equation}
From Eq.~\ref{eq2-03}, each player's payoff is obtained as
\begin{equation}
\begin{split}
	&u_1^{\mathrm{BB}}:=u_1(x_1^{\mathrm{BB}},x_2^{\mathrm{BB}})=r^2/(1+r)^2,\\
	&u_2^{\mathrm{BB}}:=u_2(x_1^{\mathrm{BB}},x_2^{\mathrm{BB}})=1/(1+r)^2.
\end{split}
\label{eq2-04}
\end{equation}

\subsection{LB equilibrium}
Player 1's L-response is given by
\begin{equation}
\begin{split}
	f_1^{\mathrm{L}}(x_2)&:=\mathrm{argmax}_{x_1}u_1(x_1,f_2(x_1))\\
	&=\mathrm{argmax}_{x_1}u_1(x_1,f_2^{\mathrm{B}}(x_1))\\
	&=\left\{\begin{array}{ll}
		r/4 & (1\le r<2) \\
		1/r & (2\le r) \\
	\end{array}\right.
\end{split}
\label{eq2-05}
\end{equation}
Then, we get the actions and payoffs at the crossing point as
\begin{equation}
\begin{split}
	(x_1^{\mathrm{LB}},x_2^{\mathrm{LB}},u_1^{\mathrm{LB}},u_2^{\mathrm{LB}})=\left\{\begin{array}{ll}
		\displaystyle\left(\frac{r}{4},\frac{r}{2}\left(1-\frac{r}{2}\right),\frac{r}{4},\left(1-\frac{r}{2}\right)^2\right) & (1\le r<2) \\
		\displaystyle\left(\frac{1}{r},0,1-\frac{1}{r},0\right) & (2\le r) \\
	\end{array}\right.
\end{split}
\label{eq2-06}
\end{equation}

By comparing LB with BB, we can prove $x_1^{\mathrm{LB}}>x_1^{\mathrm{BB}}$, $x_2^{\mathrm{LB}}<x_2^{\mathrm{BB}}$, $u_1^{\mathrm{LB}}>u_1^{\mathrm{BB}}$, and $u_2^{\mathrm{LB}}<u_2^{\mathrm{BB}}$ for all $r>1$.

\subsection{BL equilibrium}
Player 2's L-response is given by
\begin{equation}
\begin{split}
	f_2^{\mathrm{L}}(x_1)&:=\mathrm{argmax}_{x_2}u_2(f_1(x_2),x_2)\\
	&=\mathrm{argmax}_{x_2}u_2(f_1^{\mathrm{B}}(x_2),x_2)\\
	&=1/(4r).
\end{split}
\label{eq2-07}
\end{equation}
Then, we get the actions and payoffs at the crossing point as
\begin{equation}
\begin{split}
	(x_1^{\mathrm{LB}},x_2^{\mathrm{LB}},u_1^{\mathrm{LB}},u_2^{\mathrm{LB}})=\left(\frac{1}{2r}\left(1-\frac{1}{2r}\right),\frac{1}{4r},\left(1-\frac{1}{2r}\right)^2,\frac{1}{4r}\right).
\end{split}
\label{eq2-08}
\end{equation}

By comparing BL with BB, we can prove $x_1^{\mathrm{BL}}<x_1^{\mathrm{BB}}$, $x_2^{\mathrm{BL}}<x_2^{\mathrm{BB}}$, $u_1^{\mathrm{BL}}>u_1^{\mathrm{BB}}$, and $u_2^{\mathrm{BL}}>u_2^{\mathrm{BB}}$ for all $r>1$.

\subsection{Difference between the Nash and Stackelberg equilibria}
Here, Eq.~\ref{eq1-05} can predict whether the Stackelberg and Nash equilibria are equal. From the comparison between LB and BB, we calculate Eq.~\ref{eq1-05} as
\begin{equation}
\begin{split}
	&\left.\frac{\partial u_1}{\partial x_2}\right|_{\mathrm{BB}}=-1\neq 0,\\
	&\left.\frac{\partial^2 u_2}{\partial x_1\partial x_2}\right|_{\mathrm{BB}}=(r+1)(-r+1)=\left\{\begin{array}{ll}
		=0 & (r=1)\\
		\neq 0 & (r>1)\\
	\end{array}\right.
\end{split}
\label{eq2-13}
\end{equation}
This equation shows that LB coincides with BB in the case of $r=1$ and not in the case of $r>1$. We discuss the overlap between BL and BB in the same way.

\subsection{Analytical solution of the functional dynamics}
Here, we derive the analytical solution of the functional dynamics. From Eq.~\ref{eq1-09}, the first-order term $x_1^{\mathrm{eq}*}$ of player 1's function $f_1^{*}(x_2)$ is given by
\begin{equation}
\begin{split}
	&\left.\frac{\partial u_1}{\partial x_1}\right|_{\mathrm{eq}*}+\epsilon_1a_2^{*}\left.\frac{\partial u_1}{\partial x_2}\right|_{\mathrm{eq}*}=0\\
	&\Leftrightarrow x_1^{\mathrm{eq}*}=\frac{1}{r+\epsilon_1a_2^{*}}\left\{\sqrt{r(x_2^{\mathrm{eq}*}-\epsilon_1a_2^{*}x_1^{\mathrm{eq}*})}-(x_2^{\mathrm{eq}*}-\epsilon_1a_2^{*}x_1^{\mathrm{eq}*})\right\}
\end{split}
\label{eq2-09}
\end{equation}
by substituting
\begin{equation}
\begin{split}
	&\left.\frac{\partial u_1}{\partial x_1}\right|_{\mathrm{eq}*}=\frac{rx_2^{\mathrm{eq}*}}{rx_1^{\mathrm{eq}*}+x_2^{\mathrm{eq}*}}-1,\\
	&\left.\frac{\partial u_1}{\partial x_2}\right|_{\mathrm{eq}*}=-\frac{rx_1^{\mathrm{eq}*}}{rx_1^{\mathrm{eq}*}+x_2^{\mathrm{eq}*}}.\\
\end{split}
\label{eq2-10}
\end{equation}
Then, from Eq.~\ref{eq1-10}, the second-order term $a_1^{*}$ of player 1's function is given by
\begin{equation}
\begin{split}
	&a_1^{*}\left(\left.\frac{\partial^2u_1}{\partial x_1^2}\right|_{\mathrm{eq}*}+2\epsilon_1a_2^{*}\left.\frac{\partial^2u_1}{\partial x_1\partial x_2}\right|_{\mathrm{eq}*}+(\epsilon_1a_2^{*})^2\left.\frac{\partial^2u_1}{\partial x_2^2}\right|_{\mathrm{eq}*}\right)+(1-\epsilon_1)\left(\left.\frac{\partial^2u_1}{\partial x_1\partial x_2}\right|_{\mathrm{eq}*}+\epsilon_1a_2^{*}\left.\frac{\partial^2u_1}{\partial x_2^2}\right|_{\mathrm{eq}*}\right)=0.\\
	&\Leftrightarrow a_1^{*}=\frac{1-\epsilon_1}{r+\epsilon_1a_2^{*}}\left\{\frac{rx_1^{\mathrm{eq}*}+x_2^{\mathrm{eq}*}}{2(x_2^{\mathrm{eq}*}-\epsilon_1a_2^{*}x_1^{\mathrm{eq}*})}-1\right\}
\end{split}
\label{eq2-11}
\end{equation}
by substituting
\begin{equation}
\begin{split}
	&\left.\frac{\partial^2 u_1}{\partial x_1^2}\right|_{\mathrm{eq}*}=-\frac{2rx_1^{\mathrm{eq}*}}{(rx_1^{\mathrm{eq}*}+x_2^{\mathrm{eq}*})^3},\\
	&\left.\frac{\partial^2 u_1}{\partial x_1\partial x_2}\right|_{\mathrm{eq}*}=\frac{r(rx_1^{\mathrm{eq}*}-x_2^{\mathrm{eq}*})}{(rx_1^{\mathrm{eq}*}+x_2^{\mathrm{eq}*})^3},\\
	&\left.\frac{\partial^2 u_1}{\partial x_2^2}\right|_{\mathrm{eq}*}=\frac{2rx_1^{\mathrm{eq}*}}{(rx_1^{\mathrm{eq}*}+x_2^{\mathrm{eq}*})^3}.
\end{split}
\label{eq2-12}
\end{equation}

\subsection{Difference between BB and LL in the functional dynamics}
From the above solution of the functional dynamics, we derive the BB and LL equilibria. First, by substituting $\epsilon_1=\epsilon_2=0$ (BB equilibrium), we get
\begin{equation}
	(x_1^{\mathrm{BB}},x_2^{\mathrm{BB}},a_1^{\mathrm{BB}},a_2^{\mathrm{BB}})=\left(\frac{r}{(1+r)^2},\frac{r}{(1+r)^2},\frac{r-1}{2r},-\frac{r-1}{2}\right).
\label{eq2-14}
\end{equation}
Then, by substituting $\epsilon_1=\epsilon_2=1$ (LL equilibrium), we get
\begin{equation}
	(x_1^{\mathrm{LL}},x_2^{\mathrm{LL}},a_1^{\mathrm{LL}},a_2^{\mathrm{LL}})=\left(\frac{r}{(1+r)^2},\frac{r}{(1+r)^2},0,0\right).
\label{eq2-15}
\end{equation}
From these equations, we confirm $x_i^{\mathrm{BB}}=x_i^{\mathrm{LL}}$, but $a_i^{\mathrm{BB}}\neq a_i^{\mathrm{LL}}$. This finding indicates that although the same equilibrium actions are achieved in BB and LL, the (gradient of the) response function between them differs.

\section{Example 2: duopoly game}
In this game, both players' payoffs are defined by
\begin{equation}
\begin{split}
	& u_1(x_1,x_2) := x_1\max(p-x_1-x_2-c_1,0), \\
	& u_2(x_1,x_2) := x_2\max(p-x_2-x_1-c_2,0).
\end{split}
\label{eq3-01}
\end{equation}

\subsection{BB equilibrium}
Both players' B-responses are given by
\begin{equation}
\begin{split}
	&x_1^{\mathrm{B}}(x_2)=\max(1-x_2-c_1,0)/2,\\
	&x_2^{\mathrm{B}}(x_1)=\max(1-x_1-c_2,0)/2.
\end{split}
\label{eq3-02}
\end{equation}

Then, both players' actions and payoffs at the BB equilibrium are given by
\begin{equation}
	(x_1^{\mathrm{BB}},x_2^{\mathrm{BB}},u_1^{\mathrm{BB}},u_2^{\mathrm{BB}})=\left(\frac{1-2c_1+c_2}{3},\frac{1-2c_2+c_1}{3},\frac{(1-2c_2+c_1)^2}{9},\frac{(1-2c_1+c_2)^2}{9}\right).
\label{eq3-03}
\end{equation}

\subsection{LB equilibrium}
Player 1's L-response is given by
\begin{equation}
\begin{split}
	f_1^{\mathrm{L}}(x_2)&=\mathrm{argmax}_{x_1}u_1(x_1,f_2^{\mathrm{B}}(x_1))\\
	&=\left\{\begin{array}{ll}
		(1-2c_1+c_2)/2 & (1>3c_2-2c_1) \\
		1-c_2 & (1\le 3c_2-2c_1) \\
	\end{array}\right.
\end{split}
\label{eq3-04}
\end{equation}
Then, we get the actions and payoffs at the crossing point as
\begin{equation}
\begin{split}
	&(x_1^{\mathrm{LB}},x_2^{\mathrm{LB}},u_1^{\mathrm{LB}},u_2^{\mathrm{LB}})\\
	&=\left\{\begin{array}{ll}
		\displaystyle\left(\frac{1-2c_1+c_2}{2},\frac{1-3c_2+2c_1}{4},\frac{(1-2c_1+c_2)^2}{8},\frac{(1-3c_2+2c_1)^2}{16}\right) & (1>3c_2-2c_1) \\
		\left(1-c_2,0,(1-c_2)(c_2-c_1),0\right) & (2\le r) \\
	\end{array}\right.
\end{split}
\label{eq3-05}
\end{equation}

By comparing Eq.~\ref{eq3-05} with \ref{eq3-03}, we get $x_1^{\mathrm{LB}}>x_1^{\mathrm{BB}}$, $x_2^{\mathrm{LB}}<x_2^{\mathrm{BB}}$, $u_1^{\mathrm{LB}}>u_1^{\mathrm{BB}}$, and $u_2^{\mathrm{LB}}<u_2^{\mathrm{BB}}$.

\subsection{BL equilibrium}
Player 2's L-response is given by
\begin{equation}
\begin{split}
	f_2^{\mathrm{L}}(x_1)&=\mathrm{argmax}_{x_2}u_2(f_1^{\mathrm{B}}(x_2),x_2)\\
	&=(1-2c_2+c_1)/2.
\end{split}
\label{eq3-06}
\end{equation}
Then, we get the actions and payoffs at the crossing point as
\begin{equation}
\begin{split}
	(x_1^{\mathrm{BL}},x_2^{\mathrm{BL}},u_1^{\mathrm{BL}},u_2^{\mathrm{BL}})=\left(\frac{1-3c_1+2c_2}{4},\frac{1-2c_2+c_1}{2},\frac{(1-3c_1+2c_2)^2}{16},\frac{(1-2c_2+c_1)^2}{8}\right).
\end{split}
\label{eq3-07}
\end{equation}

By comparing Eq.~\ref{eq3-07} with Eq.~\ref{eq3-03}, we get $x_1^{\mathrm{BL}}<x_1^{\mathrm{BB}}$, $x_2^{\mathrm{BL}}>x_2^{\mathrm{BB}}$, $u_1^{\mathrm{BL}}<u_1^{\mathrm{BB}}$, and $u_2^{\mathrm{BL}}>u_2^{\mathrm{BB}}$.

\subsection{Difference between the Nash and Stackelberg equilibria}
We now discuss whether the Nash and Stackelberg equilibria are equal. From Eq.~\ref{eq1-05}, we get
\begin{equation}
\begin{split}
	&\left.\frac{\partial u_1}{\partial x_2}\right|_{\mathrm{BB}}=\frac{(1-2c_1+c_2)^2}{3}\neq 0,\\
	&\left.\frac{\partial^2 u_2}{\partial x_1\partial x_2}\right|_{\mathrm{BB}}=-1\neq 0.
\end{split}
\label{eq3-08}
\end{equation}
Hence, the LB and BB equilibria are not equal.

\subsection{Analytical solution of the functional dynamics}
In the duopoly game, player 1'''s payoff $u_1(x_1,x_2)$ depends on $x_1$ ($x_2$) at most quadratically (linearly), and $u_2(x_1,x_2)$ depends on $x_2$ ($x_1$) at most quadratically (linearly), respectively. Thus, we get
\begin{equation}
\begin{split}
	&f_1^{*}(x_2)=\max(a_1^{*}x_2+b_1^{*},0),\\
	&f_2^{*}(x_1)=\max(a_2^{*}x_1+b_2^{*},0).
\end{split}
\label{eq3-09}
\end{equation}
From Eqs.~\ref{eq1-09} and \ref{eq1-10}, we obtain both players' functions as
\begin{equation}
\begin{split}
	&a_1^{*} = \frac{(2-\epsilon_1+\epsilon_2)-\sqrt{(2-\epsilon_1-\epsilon_2)^2+4\epsilon_1\epsilon_2}}{4\epsilon_2}, \\
	&a_2^{*} = \frac{(2-\epsilon_2+\epsilon_1)-\sqrt{(2-\epsilon_1-\epsilon_2)^2+4\epsilon_1\epsilon_2}}{4\epsilon_1}, \\
	&b_1^{*} = \frac{2(1-\epsilon_2a_1^{*})(1-c_1)-\epsilon_1(1-c_2)}{4(1-\epsilon_1a_2^{*})(1-\epsilon_2a_1^{*})-\epsilon_1\epsilon_2}, \\
	&b_2^{*} = \frac{2(1-\epsilon_2a_2^{*})(1-c_2)-\epsilon_2(1-c_1)}{4(1-\epsilon_1a_2^{*})(1-\epsilon_2a_1^{*})-\epsilon_1\epsilon_2}.
\end{split}
\label{eq3-10}
\end{equation}

\begin{figure}[htb]
\begin{center}
\includegraphics[width=0.8\linewidth]{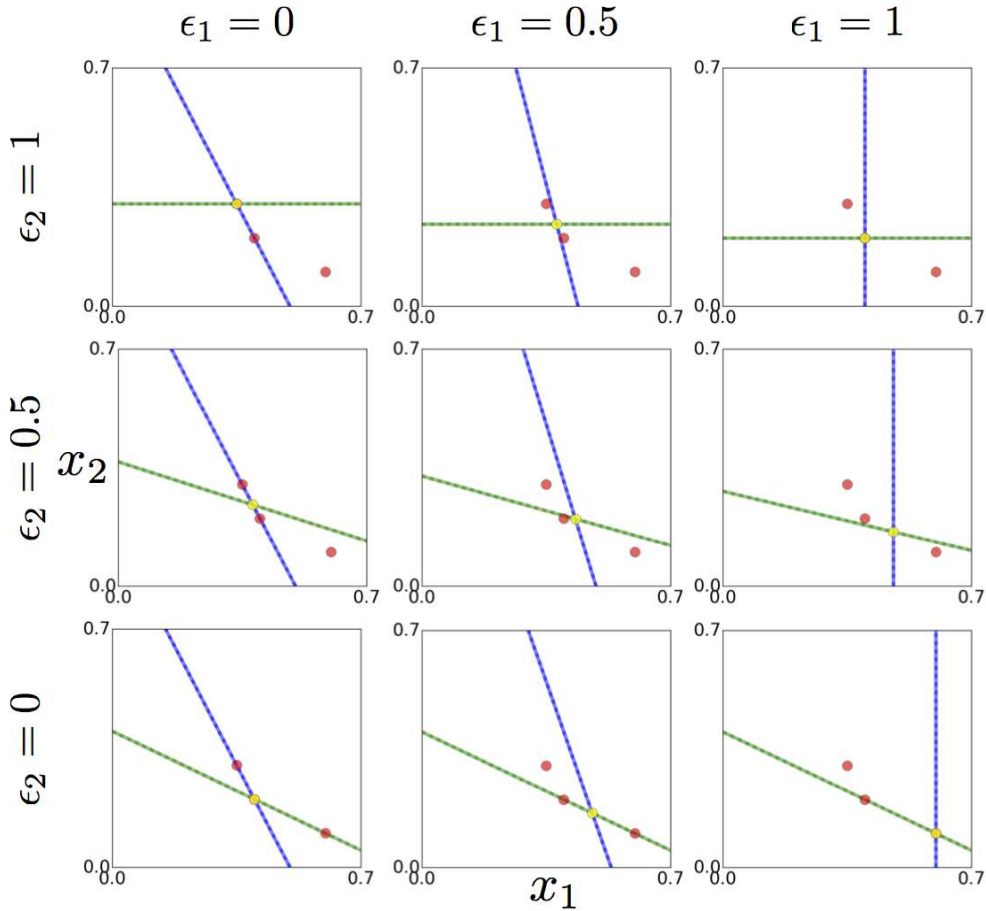}
\caption{Equilibrium functions of mutual intention recognition in duopoly game $(c_1,c_2)=(0,0.2), p=1$. For all nine figures, the X-axis (Y-axis) indicates player 1's (2's) action, denoted by $x_1$ ($x_2$). The blue (green) line indicates 1's (2's) intention and the solid (broken) line indicates the simulated (analytical) solution. These agree with each other. The yellow dot is the crossing point, while the red dots are the Nash and Stackelberg equilibria plotted for reference. The left, center, and right figures are respectively the cases of $\epsilon_1=0,0.5,1$, while the upper, center, and lower figures are respectively the cases of $\epsilon_2=0,0.5,1$.}
\label{picS01}
\end{center}
\end{figure}

\begin{figure}[htb]
\begin{center}
\includegraphics[width=0.8\linewidth]{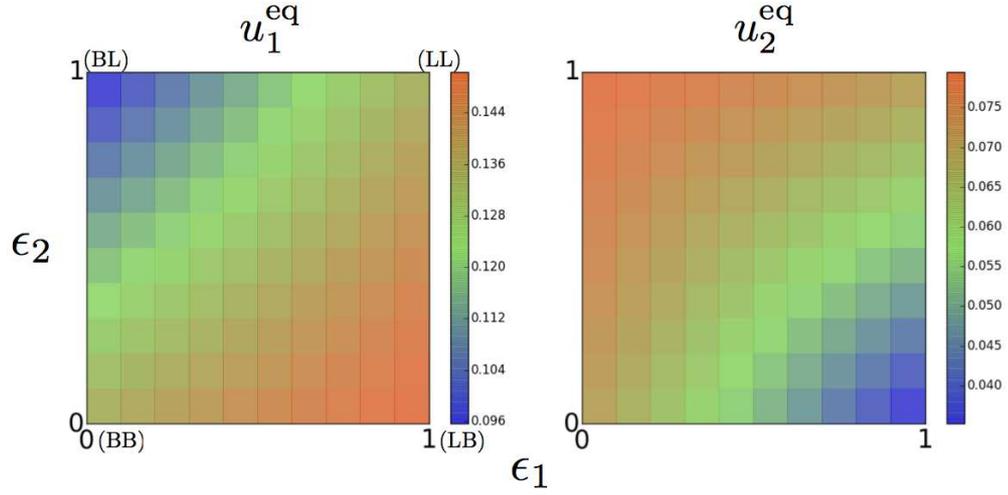}
\caption{Player 1's (left) and 2's (right) simulated payoffs in the equilibrium of intention recognition in duopoly game $(c_1,c_2)=(0,0.1)$. For all figures, the X-axis (Y-axis) indicates player 1's (2's) recognition degree.}
\label{picS02}
\end{center}
\end{figure}

Fig.~S\ref{picS01} shows that the equilibrium functions $f_1^{*}(x_2),f_2^{*}(x_1)$ obtained from the simulation agree well with the above analytic expectation. Fig.~S\ref{picS02} shows the equilibrium payoffs achieved by both players' functions $f_1^{*}(x_2),f_2^{*}(x_1)$ for $(\epsilon_1,\epsilon_2)$ with $\epsilon_1,\epsilon_2=0,0.5,1$.

% Figure S01

% Figure S02

\subsection{Simulation of the learning process}
\begin{figure}[htb]
\begin{center}
\includegraphics[width=0.4\linewidth]{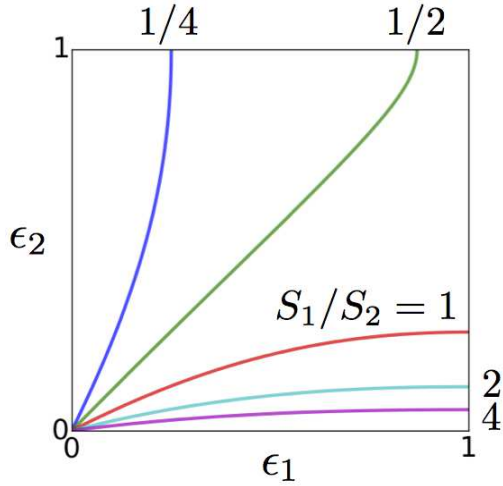}
\caption{Trajectories of the set of recognition degrees $(\epsilon_1, \epsilon_2)$ for diverse sets of learning speed $(S_1,S_2)$ in duopoly game $(c_1,c_2)=(0,0.1), p=1$. The blue, green, red, cyan, and magenta lines represent the cases of $S_1/S_2=0.25, 0.5, 1, 2, 4$, respectively. For all trajectories, the learning dynamics start from $(\epsilon_1,\epsilon_2)=(0,0)$.
}
\label{picS03}
\end{center}
\end{figure}

\begin{figure}[htb]
\begin{center}
\includegraphics[width=0.4\linewidth]{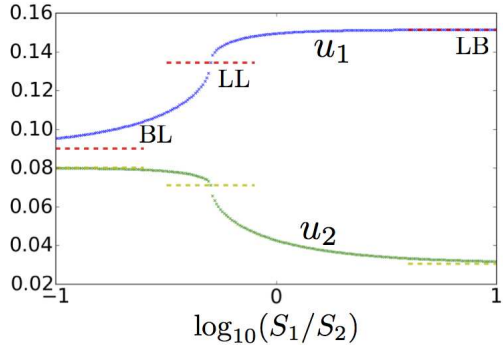}
\caption{Player 1's payoff (blue dots) and player 2's payoff (green dots) in the equilibrium state as a function of $S_1/S_2$, i.e., for the ratio of learning speeds $(S_1,S_2)$ in duopoly game $(c_1,c_2)=(0,0.1), p=1$. The red (yellow) broken lines indicate player 1's (2's) payoff at the LB (right), LL (middle), and BL (left) equilibria for reference, respectively.
}
\label{picS04}
\end{center}
\end{figure}

Here, we consider a process to change the degree of recognition $\epsilon_i$ to increase the first agent's payoff under the other's intention. Fig.~S\ref{picS03} shows how both players' recognition degrees change for some $S_1/S_2$. In addition, the finally achieved payoffs are plotted in Fig.~S\ref{picS04}.

% Figure S03

% Figure S04

\section{Example 3: probabilistic prisoner's dilemma game}
Our third example is a prisoner's dilemma game in which each of the players chooses whether to cooperate or defect. In general, each player's payoff is given by
\begin{equation}
\begin{split}
	& u_1(x_1,x_2)=T(1-x_1)x_2+Rx_1x_2+P(1-x_1)(1-x_2)+Sx_1(1-x_2),\\
	& u_2(x_1,x_2)=T(1-x_2)x_1+Rx_2x_1+P(1-x_2)(1-x_1)+Sx_2(1-x_1).
\end{split}
\label{eq4-01}
\end{equation}
Here, both $T>R>P>S$ and $2R>T+S$ are required. Then, we consider a situation that each player determines the probability $x_1,x_2$ to cooperate as his/her action.

First, each player's B-response is given by
\begin{equation}
\begin{split}
	& f_1^{\mathrm{B}}(x_2)=0, \\
	& f_2^{\mathrm{B}}(x_1)=0.
\end{split}
\label{eq4-02}
\end{equation}
Eq.~\ref{eq4-02} indicates that every player is better off choosing defection regardless of the opponent's cooperativeness. In other words, both B-responses are constant compared with the change in the opponent's action. Then, we get the following actions and payoffs at the BB equilibrium:
\begin{equation}
	(x_1^{\mathrm{BB}},x_2^{\mathrm{BB}},u_1^{\mathrm{BB}},u_2^{\mathrm{BB}})=(0,0,P,P).
\label{eq4-03}
\end{equation}

Player 1's L-response is given by
\begin{equation}
\begin{split}
	f_1^{\mathrm{L}}(x_2)&=\mathrm{argmax}_{x_1}u_1(x_1,f_2^{\mathrm{B}}(x_1))\\
	&=0.
\end{split}
\label{eq4-04}
\end{equation}
This equation indicates that player 1's L-response is equal to the B-response. Hence, player 1 gains no advantage by one-way recognition. This is also true for player 2's one-way recognition. Hence, the LB, BL, and BB equilibria are all given by $x_1^{\mathrm{eq}}=x_2^{\mathrm{eq}}=0$, and all of them agree.

\end{document}